\documentclass[12pt]{article}
\usepackage[T2A]{fontenc}
\usepackage[cp1251]{inputenc}
\usepackage{amsfonts}
\usepackage{eufrak}
\usepackage{amssymb}
\usepackage{amsmath}


\begin{document}

\centerline{\bf On a functional equation appearing in characterization of} 
\centerline{\bf distributions by the optimality of  an estimate\footnote
{This research was supported in part by the grant “Network of Mathematical
        Research 2013--2015”}}

\vskip 1 cm

\centerline{G.M. FELDMAN}

\bigskip

\centerline{B. Verkin Institute for Low Temperature Physics  and Engineering}
\centerline{ Kharkiv,
 Ukraine} 

\bigskip

\centerline{and}

\bigskip

\centerline{P. Graczyk}

\bigskip

\centerline{  Laboratoire LAREMA UMR 6093 - CNRS}
\centerline{Universite a'Angers, France} 

\vskip  1 cm

  \makebox[20mm]{ }\parbox{125mm}{ \small
Let  $X$ be a second countable 
 locally compact Abelian group containing no subgroup topologically isomorphic 
 to the circle group
$\mathbb{T}$. Let $\mu$ be a probability distribution on $X$
such that its characteristic function  $\widehat\mu(y)$ does not vanish and 
  $\widehat\mu(y)$ for some  $n \geq 3$ satisfies the equation
  $$
  \prod_{j=1}^{n} \hat\mu(y_j + y) = \prod_{j=1}^{n}
\hat\mu(y_j - y), \quad \sum_{j=1}^{n} y_j = 0, \quad
y_1,\dots,y_n,y \in Y.
$$ 
    Then $\mu$ is a convolution of a Gaussian distribution and a distribution
    supported in the subgroup of $X$ generated by elements of order 2.}

\vskip  1 cm

The present note is devoted to study a functional equation on a locally compact Abelian group
group which appears in characterization of probability  distributions by the optimality of  
an estimate.

Let  $X$ be a second countable 
 locally compact Abelian group, $Y=X^*$ be its character group, $(x,y)$
be the value of a character  $y\in Y$ at an element  $x\in X$.
Denote by ${\rm M}^1(X)$ the convolution semigroup of probability distribution on
the group  $X$, and denote by $$\widehat\mu(y) =
\int_{X}(x, y)d \mu(x)$$ the characteristic function of a distribution  $\mu\in {\rm M}^1(X)$.
For $\mu \in {\rm M}^1(X)$
define the distribution  $\bar \mu \in {\rm M}^1(X)$ by the
formula $\bar \mu(B) = \mu(-B)$ for all Borel sets of $X$. Then $\widehat{\bar\mu}(y)=\overline{\widehat\mu(y)}.$

A distribution  $\gamma\in {\rm M}^1(X)$ is called Gaussian if its characteristic function is represented
in the form
\begin{equation}\label{0}
\widehat\gamma(y)= (x,y)\exp\{-\varphi(y)\},
\end{equation}
where $x \in X$, and $\varphi(y)$ is a continuous nonnegative function on the group $Y$ 
satisfying the equation
\begin{equation}\label{1}
\varphi(y_1 + y_2) + \varphi(y_1
- y_2) = 2[\varphi(y_1) + \varphi(y_2)], \quad y_1, \ y_2 \in
Y. 
\end{equation}
A Gaussian distribution is called symmetric if in (\ref{0}) $x=0$.
Denote by $\Gamma(X)$ the set of Gaussian distributions on the group $X$.

Consider a probability space $(X, {\cal B}, \mu)$, where
${\cal B}$ is a $\sigma$-algebra of Borel subsets of $X$, and
$\mu \in {\rm M}^1(X)$.
Form a family of distributions 
$\mu_\theta(A) = \mu(A - \theta), \ A \in {\cal B}, \ \theta
\in X$.  Denote by $\Pi$ a class of estimates 
$f : X^n \mapsto X$ satisfying the condition
$f(x_1 + c, \dots , x_n + c) = f(x_1, \dots, x_n) + nc$
for all $x_1,\dots ,x_n, c \in X$.
According to  \cite{Ru1} (see also \cite{Ru2},  \cite[Ch. 7, \S 7.10]{KaLiRa}),
an estimate $f_0\in\Pi$ of a parameter $n\theta$ is called an optimal estimate in the class $\Pi$ 
for a sample volume $n$ if for any estimate  $f\in\Pi$ and for all $y \in Y$ the inequality
$$ {\bf E}_\theta|(f_0({\bf x}), y) - (n\theta,
y)|^2 \leq {\bf E}_\theta|(f({\bf x}), y) - (n\theta, y)|^2 
$$
holds.
It turns out that the existence of an optimal estimate of  the papameter 
$n\theta$ gives the possibility in some cases to describe completely the possible distributions $\mu$.

As has been proved in (\cite{Ru1}), if an estimate $f_0$ is represented in the form
\begin{equation}\label{2}
(f_0({\bf x}), y) = (f({\bf x}), y) \frac{{\bf E}_0[(f({\bf x}), -y)| {\bf z}]}
{|{\bf E}_0[(f({\bf x}), -y)| {\bf z}]|}, \quad y \in Y, 
\end{equation}
where $f\in \Pi$, and ${\bf z} = (x_2 - x_1,\dots ,x_n - x_1)$, then
$f_0\in \Pi$,   $f_0$ does not depend on the choice of
$f$ and $f_0$ is an optimal estimate.
It follows from (\ref{2}) that 
$f_0$ is an optimal  estimate if and only if 
$\arg{\bf E}_0[(f_0({\bf x}), y)| {\bf z}] = 0$ .
When $f_0({\bf x}) = \sum_{j=1}^{n}x_j$
 it follows from this that the characteristic function $\widehat\mu(y)$ satisfies
 the equation
\begin{equation}\label{3}
\prod_{j=1}^{n} \hat\mu(y_j + y) = \prod_{j=1}^{n}
\hat\mu(y_j - y), \quad \sum_{j=1}^{n} y_j = 0, \quad
y_1,\dots,y_n,y \in Y,
\end{equation}
and $\widehat\mu^{n}(y) > 0$. When $\ n \geq 3 $
this implies that if a group  $X$ contains no elements of order 2, then
 $\mu \in \Gamma(X)$ (see \cite{Ru1}).

This note is devoted to solving of equation (\ref{3}) in a general case when $X$ is a locally
compact Abelian group. Let us fix the notation.
Denote by $f_2 : X \mapsto X$  the endomorphism of $X$ defined by the formula 
 $f_2(x)=2x$.  Put $X_{(2)} = {\rm Ker} \
f_2,$ $X^{(2)} = {\rm Im} \ f_2.$ Denote by $\mathbb{T}$
the circle group, and by  $\mathbb{Z}$ the group of integers.

Let $\psi(y)$ be an arbitrary function on the group  $Y$ and $h \in
Y$.  Denote by $\Delta_h$ the finite difference operator
$$
\Delta_h \psi(y) = \psi(y + h) - \psi(y), \quad y \in Y.
$$
A continuous function $\psi(y)$ on the group   $Y$ is called a polynomial 
if
\begin{equation}\label{4}
\Delta_h^{m+1} \psi(y) = 0 
\end{equation}
for some $m$ and for all  $y, h \in Y$.  The minimal $m$ for which (\ref{4}) holds
is called the degree of the polynomial  $\psi(y)$.

From analytical point of view the result proved in \cite{Ru1}
can be reformulated in the following way. Let
$\mu\in {\rm M}^1(X)$, the characteristic function 
 $\widehat\mu(y)$ satisfy equation (\ref{3}) for some 
 $n \geq 3$ and $\widehat\mu^n(y)>0$.  Then if the group $X$
contains no elements of order 2, then  $\mu \in \Gamma(X)$.

It is easy to see that if $\gamma$ is a symmetric Gaussian distribution on the group $X$
and $\pi \in {\rm M}^1(X_{(2)})$, 
then the characteristic functions $\widehat\gamma(y)$ and $\widehat\pi(y)$
satisfy equation (\ref{3}), and hence  the characteristic function of
the distribution  $\mu = \gamma*\pi$ also satisfies  equation (\ref{3}).
Describe first the groups $X$ for which the converse statement is true.

\bigskip

{\bf Theorem 1}. {\it Let  $X$ be a second countable 
 locally compact Abelian group, $\mu \in {\rm M}^1(X)$.
  Let the characteristic function  $\widehat\mu(y)$ satisfy equation $(\ref{3})$ 
  for some  $n \geq 3$ and $\widehat\mu(y) \ne 0$. 
  Assume that the following condition holds:  $(i)$
the group $X$ contains no subgroup topologically isomorphic to the circle group
$\mathbb{T}$. Then $\mu=\gamma*\pi$, where $\gamma\in \Gamma(X)$ and $\pi\in{\rm M}^1(X_{(2)})$}.

\bigskip

{\bf Proof}.  Set $\nu = \mu*\bar\mu$. Then  $\widehat\nu(y) =
|\widehat\mu(y)|^2 > 0$.  Put $\psi(y) = - \ln\widehat\nu(y)$.
Equation (\ref{3}) is equivalent to the equation
\begin{equation}\label{5}
\sum_{j=1}^{n}\psi(y_j + y) =  \sum_{j=1}^{n}\psi(y_j - y),
\quad  \sum_{j=1}^{n}y_j = 0, \quad y_1,\dots,y_n, y \in Y. 
\end{equation} 	
We also note that 
\begin{equation}\label{6}
\psi(-y) = \psi(y), \quad y \in Y.
\end{equation}  
Substituting in (\ref{5}) $y_3=-y_1-y_2, \quad y_4=\dots=y_n=0$
and taking into account (\ref{6}), we get 
\begin{equation}\label{7} \psi(y_1 + y_2 + y) - \psi(y_1 +
y_2 - y) = \psi(y_1 + y) - \psi(y_1 - y) $$ $$ +  \psi(y_2 + y)
 - \psi(y_2 - y), \quad y_1, y_2, \ y \in Y.  
\end{equation} 
Setting successively   $y = y_1 + y_2$, $y = y_1$, $y =
y_2$, we find from (\ref{7}) that 
$$ 
\psi(2y_1 + 2y_2) = \psi(2y_1) + 2\psi(y_1 +
y_2) - 2\psi(y_1 - y_2) + \psi(2y_2), \quad y_1, y_2 \in Y.
$$
This implies that 
$$
\psi(2y_1 + 2y_2) + \psi(2y_1 - 2y_2) =
2[\psi(2y_1) + \psi(2y_2)], \quad y_1, y_2 \in Y,
$$
i.e. the function 
 $\psi(y)$ satisfies equation (\ref{1}) on the subgroup  $Y^{(2)}$, 
 and hence, the function 
 $\psi(y)$ satisfies equation (\ref{1}) on the subgroup  $\overline{Y^{(2)}}$.
  Denote by  $\varphi_0(y)$ the restriction of the function  $\psi(y)$
to the subgroup $\overline{Y^{(2)}}$.

It is well known that we can associate to each function $\varphi(y)$ satisfying equation 
(\ref{1}) a symmetric 2-additive function
$$
\Phi(u, v) = \frac{1}{2} [ \varphi(u + v) - \varphi(u) -
\varphi(v)], \quad u, v \in Y.
$$
Then $\varphi(y) = \Phi(y, y)$. Using this representation it is not difficult to verify that
the  function 
$\psi(y)$ on the subgroup  $\overline{Y^{(2)}}$ satisfies the equation
\begin{equation}\label{8}
\Delta_k^2\Delta_{2h} \ \psi(y) = 0, \quad k, \ y \in
\overline{Y^{(2)}}, \quad h \in Y.  
\end{equation}

Return to equation (\ref{7}) and apply the finite difference method 
to solve it.  Let $h_1$  be an arbitrary element of the group  $Y$. 
Put $k_1 = h_1$. Substitute  $y_2+h_1$ for $y_2$ and $y+k_1$ for $y$
in equation (\ref{7}).  Subtracting equation  (\ref{7}) from the resulting equation 
we obtain 
\begin{equation}\label{9}
\Delta_{2h_1} \ \psi(y_1 + y_2 + y ) = \Delta_{h_1} \
\psi(y_1 + y) - \Delta_{-h_1} \ \psi(y_1 - y) +
\Delta_{2h_1} \ \psi(y_2 + y). 
\end{equation}
Next, let $h_2$ be an arbitrary element of the group $Y$. Put $k_2 = - h_2$. 
Substitute  $y_2+h_2$ for $y_2$ and $y+k_2$ for $y$
in equation (\ref{9}). Subtracting equation  (\ref{9}) from the resulting equation 
we get
$$
\Delta_{-h_2}\Delta_{h_1} \ \psi(y_1 + y ) -
 \Delta_{h_2}\Delta_{-h_1} \ \psi(y_1 - y ) = 0.
$$
Reasoning similarly we find from this 
$$
\Delta_{2h_3}\Delta_{-h_2}\Delta_{h_1} \ \psi(y_1 + y ) = 0,
$$
and finally 
\begin{equation}\label{a1}
\Delta_{2h_3}\Delta_{-h_2}\Delta_{h_1}\Delta_{y_1} \
\psi(y) = 0.
\end{equation}
Note that $h_j,  y_1,   y$ are arbitrary elements of the group $Y$.
Setting in  (\ref{a1}) $h_3=h$,  $-h_2=h_1=y_1=k$,
we find
\begin{equation}\label{10}
\Delta_k^3\Delta_{2h} \ \psi(y) = 0, \quad k, \ h, \ y \in
 Y.  
\end{equation}
Fix  $h \in Y$.  On  the one hand, it follows from (\ref{10}) that the function 
$\Delta_{2h}\psi(y)$ is a polynomial of degree  $\leq 2$ on the group
$Y$.  On the other hand, as follows from (\ref{8}) the function $\Delta_{2h}\psi(y)$
is a polynomial of  degree $\leq 1$ on the subgroup $\overline{Y^{(2)}}$.
 Then as not difficult to verify, the function  
$\Delta_{2h}\psi(y)$ must be a polynomial of  degree $\leq 1$ on the group  $Y$, i.e.
\begin{equation}\label{11}
\Delta_k^2\Delta_{2h} \ \psi(y) = 0, \quad k, \ h, \ y \in Y.
\end{equation}

Theorem 1 follows now from the following lemma.

\bigskip

{\bf Lemma 1} (\cite[Prop. 1]{Fe2}). {\it Let  $X$ be a second countable 
 locally compact Abelian group containing no subgroup topologically isopmorphic to 
 the circle group $\mathbb{T}$. Let $\mu \in {\rm M}^1(X)$, $\nu=\mu*\bar\mu$ and
$$ \hat\nu(y)=\exp\{-\psi(y)\}, $$
where the function 
$\psi(y)$ satisfies equation  $(\ref{11})$.
Then $\mu=\gamma*\pi$, where $\gamma\in \Gamma(X)$ and $\pi\in{\rm M}^1(X_{(2)})$}.

\bigskip

 {\bf Remark 1}. Obviously, the above mentioned Rukhin's theorem follows directly from Theorem 1.

 \bigskip

  {\bf Remark 2}. Let  $X$ be a second countable 
 locally compact Abelian group containing a subgroup topologically isomorphic to 
 the circle group $\mathbb{T}$. Then we can consider any distribution $\mu$ on the circle group 
 $\mathbb{T}$  as 
 a distribution on  $X$. Note that  $\mathbb{Z}$ is the character group of $\mathbb{T}$.
 Following to \cite{Ru1} consider on the group  $\mathbb{Z}$ the function
$$f(m)=\begin{cases}
\exp\{-m^2\}, & \text{\ if\ }\ \ m \in
\mathbb{Z}^{(2)},
\\ \exp\{-m^2 + \varepsilon\}, & \text{\ if\ }\ m \not\in
\mathbb{Z}^{(2)},
\end{cases}
$$
where  $\varepsilon > 0$ is small enough.  Then
$$ \rho(t) = \sum_{m=-\infty}^\infty f(m)e^{-imt} > 0.$$ 
Let $\mu$ be a distribution on $\mathbb{T}$ with density $ \rho(t)$ with respect to the Lebesque
measure. Then $f(m)$ is the characteristic function of a distribution
  $\mu$ on the circle group $\mathbb{T}$.  Considering $\mu$ as a distribution on the group
  $X$, we see that  $\widehat\mu(y) >0$ and the characteristic function  $\widehat\mu(y)$ 
  satisfies equation (\ref{3}), but as easily seen,
  $\mu\not\in \Gamma(X)*{\rm M}^1(X_{(2)})$.  This example shows that condition  $(i)$ 
  in Theorem 1 is sharp.

\bigskip

{\bf Remark 3}. Let  $X$ be a second countable 
 locally compact Abelian group. In the articles \cite{Fe1} and \cite{Fe2} 
 (see also \cite[\S 16]{Fe2}) were studied group analogs of 
 the well-known Heyde theorem, where a Gaussian distribution is characterized 
 by the symmetry of the conditional distribution of a linear form 
 $L_2 =\beta_1\xi_1 + \cdots + \beta_n\xi_n$ of independent random variables $\xi_j$
 given  $L_1 =\alpha_1\xi_1 + \cdots + \alpha_n\xi_n$ (coefficients of the forms are 
 topological automorphisms of the group  $X$).
Let $\widehat\mu_j(y)$ be the characteristic function of the random variable 
 $\xi_j.$ It is interesting to remark that if the number of 
 independent random variables $n=2$, then the functions  
$\psi_j(y)=-\ln|\widehat\mu_j(y)|^2$ also satisfy equation (\ref{11}).
For the groups $X$ containing no subgroup topologically isomorphic to the circle group
 $\mathbb{T}$, and also for the two-dimensional torus $X=\mathbb{T}^2$ 
 this implies that all $\mu_j \in
 \Gamma(X)*{\rm M}^1(X_{(2)})$.

\bigskip

We use Theorem 1 to prove the following statement, a  significant part of which 
refers to the case when the group  $X$ contains a subgroup topologically isomorphic to
the circle group  $\mathbb{T}$.

\bigskip

{\bf Theorem 2}. {\it Let  $X$ be a second countable 
 locally compact Abelian group. Let $\mu \in {\rm M}^1(X)$,
let the characteristic function $\widehat\mu(y)$ satify equation $(\ref{3})$ for some odd
$n$, and $\widehat\mu^n(y) > 0$.  Assume that the group $X$ satisfies the condition: $(i)$
the subgroup $X_{(2)}$ is finite. Then $\mu = \gamma_0*\pi$, where
$\gamma_0 \in \Gamma(X)$, and $\pi$ is a signed measure on  $X_{(2)}$.}

\bigskip

{\bf Proof.}  Put $\psi(y) = -\ln|\widehat\mu(y)|$. Then the function 
$\psi(y)$ satisfies equation (\ref{11}). As has been proved in \cite{Fe2} in this case 
the function  $\psi(y)$ is represented in the form 
$$ \psi(y) = \varphi(y) + r_\alpha, \quad y \in
y_\alpha + \overline{Y^{(2)}},
$$
where $\varphi(y)$ is a continuous  function satisfying equation (\ref{1}), and 
 $Y =
\bigcup_\alpha {(y_\alpha + \overline{Y^{(2)}})}$ is a decomposition of the group $Y$ with respect to
the subgroup
 $\overline{Y^{(2)}}$.
  Since $X_{(2)}$ is a finite subgroup, it is easy to see that the function
 $g(y)
=\exp\{ - r_\alpha\}, \quad y \in
y_\alpha + \overline{Y^{(2)}},$
    is the chracteristic function of a signed measure  $\pi$ on the subgroup $X_{(2)}$. 
    It follows from this that
   $$
   |\hat\mu(y)| = \hat\gamma(y) \hat\pi(y),
   $$
where $\gamma \in \Gamma(X)$ and $\widehat\gamma (y)=
\exp\{-\varphi(y)\}$.

   Set $l(y) = \widehat\mu(y) / |\hat\mu(y)|$ and check that the function
 $l(y)$ is a character of the group $Y$.  Hence, Theorem 2 will be proved.

  Note that the function $l(y)$ satisfies equation (\ref{3}) and 
 \begin{equation}\label{12}
   l(-y) = \overline {l(y)}, \quad l^n(y) = 1, \quad y \in Y.
  \end{equation}
  Put in (\ref{3}) $y_2 = -y_1$, $y_3 =\dots= y_n = 0$. We get
   $$
   l^{n-2}(y) l(y_1 + y) l(-y_1 + y) = l^{n-2}(-y) l(y_1 -
   y) l(-y_1 - y), \quad y, \ y_1, \ y_2 \in Y.
   $$
   Taking into account (\ref{12}), it follows from this that
   $$
   l^{2}(y + y_1) l^{2}(y - y_1) = l^{4}(y),  \quad y, \ y_1 \in
Y.
$$
Set $m(y) = l^2(y)$.  Then the function  $m(y)$
satisfies the equation
\begin{equation}\label{13}
m(u + v) m(u - v) = m^2(u), \quad u, \ v \in Y. 
\end{equation}
We find by induction from (\ref{13}) that
\begin{equation}\label{14}
m(py) = m^p(y), \quad p \in \mathbb{Z}, \ y \in Y.  
\end{equation}

Now we formulate as a lemma the following statement.

\bigskip

{\bf Lemma 2.} {\it Let $Y = Y_1 + Y_2$, let a continuous function  $m(y)$ 
on $Y$ satisfy equation $(\ref{13})$ and
$m^{n}(y) = 1$ for some odd $n$. Then, if the restriction of the function 
$m(y)$ to $Y_j$ is a character of the group $Y_j$, $j =
1, 2$, then $m(y)$ is a character of the group $Y$.}

\bigskip

{\bf Proof}.  Denote by  $y = (y_1, y_2)$, $y_1\in Y_1, y_2\in Y_2$ elements of the group $Y$.  
Put $a(y_1, y_2) = m(y_1, 0) m(0, y_2)$, $b(y_1,
y_2) = m(y_1, y_2)/a(y_1, y_2)$.  Then $b(y_1, 0) = b(0, y_2)
= 1$, $y_1\in Y_1, y_2\in Y_2$.  It is obvious that the function  
$b(y_1, y_2)$ also satisfies equation 
  (\ref{13}).  Substitute in (\ref{13}) $u = (y_1, 0)$, $v = (y_1,
y_2)$.  We have
$$
b(2y_1, y_2) b(0, -y_2) = b^2(y_1, 0), \quad y_1\in Y_1, y_2\in Y_2.
$$
This implies that
$b(2y_1, y_2) = 1$ for $y_1\in Y_1, y_2\in Y_2$.  In particular, 
$b(2y_1, 2y_2) = 1$.  But it follows from (\ref{14}) that $b(2y_1, 2y_2) =
b^2(y_1, y_2)$.  Hence, $b(y_1, y_2) = \pm 1$.
Since $b^n(y_1, y_2) = 1$ and $n$ is odd, we have $b(y_1,
y_2)=1$ for $y_1\in Y_1, y_2\in Y_2$, i.e.  $m(y_1, y_2) = a(y_1, y_2)$ is a character of the group $Y$.

\bigskip

Continue the proof of Theorem 2. Since, by the assumption,  
$X_{(2)}$ is a finite subgroup, there exist $q \ge 0$ such that the group  $X$
contains a subgroup topologically isomorphic to the group $\mathbb{T}^q$, but $X$ does not  
contain  a subgroup topologically isomorphic to the group $\mathbb{T}^{q+1}$.
It is well known that a subgroup of $X$ topologically isomorphic to a group of the form 
$\mathbb{T}^k$ is a topologically
direct summand in $X$.  For this reason the group  $X$ is represented in the form $X =
\mathbb{T}^q + G$, where the group  $G$ contains no
subgroup topologically isomorphic to the circle group  $\mathbb{T}$. We have $Y \cong
\mathbb{Z}^q + H$, $H = G^*$.  It follows from Lemma 2 and (\ref{14}) by induction that
the function  $m(y)$ on the group  $\mathbb{Z}^q$,
satisfying equation  (\ref{13}) and the condition  $m^n(y) = 1$ is a character of
the group $\mathbb{Z}^q$. By Theorem 1 the  restriction of the function 
$\widehat\mu(y)$ to $H$ is a product of the characteristic function of a
Gaussian distribution on the group  $G$ and the characteristic function of a distribution
on the subgroup  $G_{(2)}$. Taking into account that the characteristic function 
of any distribution on $G_{(2)}$ takes only real values, it follows from the equality 
 $$\widehat\mu^2(y) = |\widehat\mu(y)|^2 m(y), \quad y\in Y,$$ that the restriction 
of the function $m(y)$ to $H$ is a character of the subgroup $H$.
Applying again Lemma 2 to the group $Y$, we obtain that $m(y)$ is a character of the group $Y$.

Since  $n$ is odd, we have $2 r + n s = 1$ for some integers $r$ and
$s$.  Taking into account (\ref{12}) this implies that $l(y) = (l(y))^{2r + ns} = (m(y))^r$ is a character 
of the group  $Y$.  Theorem 2 is completely proved.

We note that the example given in Remark 2 shows that a signed measure 
 $\pi$ needs not be a measure.

 \bigskip

{\bf Remark 4}. Consider the infinite-dimensional torus
$X=\mathbb{T}^{\aleph_0}$. Then $Y \cong {\mathbb{Z}}^{{\aleph_0}*}$,
where ${\mathbb{Z}}^{{\aleph_0}*}$ is the group of all sequences of integers
such that in each sequence only finite number of members are not equal to zero.

Consider on the group $\mathbb{Z}$ the sequence of the functions 
$$f_k(m) =\begin{cases}
 \exp\{-a_k m^2\}, & \text{\ if\ }\ \ m \in
\mathbb{Z}^{(2)},
\\ \exp\{-a_k m^2 + k\}, & \text{\ if\ }\ m \not\in
\mathbb{Z}^{(2)},
\end{cases}
$$
where $k = 1, 2, \dots$. Put
$$
f(m_1, \dots, m_l, 0, \dots) = \prod _{k=1}^{l}f_k(m_k), \quad (m_1,
\dots, m_l, 0, \dots) \in {\mathbb{Z}}^{{\aleph_0}*}.
$$
Take $a_k > 0$ such that
$$
\sum_{(m_1,\dots, m_l, 0, \dots) \in {\mathbb{Z}}^{{\aleph_0}*}}f(m_1, \dots,
m_l, 0, \dots) < 2.
$$
Then 
$$
\rho(t_1, \dots, t_l, \dots)=\sum_{(m_1,\dots, m_l, 0, \dots) \in {\mathbb{Z}}^{{\aleph_0}*}}f(m_1, \dots,
m_l, 0, \dots)e^{-i(m_1t_1+\dots+m_lt_l+\dots)}>0, \quad t_j\in {\mathbb{R}}.$$

It follows from this  that  $f(m_1, \dots, m_l, 0, \dots)$ is the characteristic function of a  
distribution $\mu \in {\rm
M}^1(\mathbb{T}^{\aleph_0})$ such that $\widehat\mu(y) >0$ and
$\widehat\mu(y)$ satisfies equation (\ref{3}), but $\mu$ can not be represented as a convolution 
 $\mu = \gamma * \pi$, where $\gamma \in
\Gamma(\mathbb{T}^{\aleph_0})$, and  $\pi$ a signed measure on the group
$\mathbb{T}^{\aleph_0}_{(2)}$. The subgroup $\mathbb{T}^{\aleph_0}_{(2)}$ is infinite in 
this case. This example shows that condition  $(i)$
in Theorem 2 is sharp.

\bigskip

{\bf Remark 5}. We assumed in Theorem 2 that 
 $n$ is odd. This condition can not be omitted  even for the circle group 
 $X = \mathbb{T}$. Indeed, let 
$n = 4$. Take $a$ in such a way that the function 
$$ f(m) = \exp\{- a m^2 + i \frac{\pi}{2} m^3\}, \quad m \in \mathbb{Z} $$ be the characteristic function of a distribution  $\mu \in
{\rm M}^1(\mathbb{T})$. On the one hand, it is ibvious that the function  $f(m)$
satisfies equation (\ref{3}) and  $f^4(m)
> 0$, $m \in \mathbb{Z}$.  On the other hand, the distribution  $\mu$ can not be represented in the form
$\mu = \gamma * \pi$, where $\gamma \in
\Gamma(\mathbb{T})$, and $\pi$ is a signed measure on $\mathbb{T}_{(2)}$.

This example also shows that a function  $f(y)$ 
satisfying equation (\ref{3}), generally speaking, needs not be real.

\end{document}